\documentstyle[amstex,amssymb,epsfig]{amsart}

\newcommand{\includeps}[1]{\psfig{#1}}

\theoremstyle{definition}
\newtheorem{thm}{Theorem}
\newtheorem{prop}{Proposition}[section]
\newtheorem{lem}[prop]{Lemma}
\newtheorem{cor}[prop]{Corollary}
\newcommand{\EE}{{\cal E}}
\newcommand{\D}{{\cal D}}
\newcommand{\F}{{\cal F}}
\newcommand{\R}{{\Bbb R}}
\newcommand{\k}{\kappa}
\newcommand{\al}{\alpha}
\newcommand{\be}{\beta}
\newcommand{\w}{\omega}

\newcommand{\di}{\partial}
\newcommand{\restr}{\negthickspace \mid}
\newcommand{\sn}{\operatorname{sn}}
\newcommand{\cn}{\operatorname{cn}}
\newcommand{\dn}{\operatorname{dn}}
\newcommand{\sign}{\operatorname{sgn}}
\newcommand{\bel}[1]{\begin{equation}\label{#1}}
\begin{document}
\title{Knot types, homotopies and stability of closed elastic rods}
\author{Thomas A. Ivey}
\address[Ivey]{Dept. of Mathematical Sciences, Ball State University,
Muncie IN 47306}
\email{tivey@@math.bsu.edu}
\author{David A. Singer}
\address[Singer]{Dept. of Mathematics, Case Western Reserve University,
Cleveland, OH 44106-7058}
\email{das5@@po.cwru.edu}
\date{Oct. 16, 1997}
%% ABSTRACT OMITTED IN PUBLISHED PAPER
\begin{abstract}
The energy minimization problem associated to uniform, isotropic, linearly 
elastic rods leads to a geometric variational problem for the rod centerline,
whose solutions include closed, knotted curves. We give a complete 
description of the space of closed and quasiperiodic solutions.
The quasiperiodic curves are parametrized by a two-dimensional disc. The
closed curves arise as a countable collection of one-parameter families, 
connecting the $m$-fold covered circle to the $n$-fold covered
circle for any $m,n$ relatively prime. Each family contains exactly one
self-intersecting curve, one elastic curve, and one closed curve of 
constant torsion.  Two torus knot types are represented in each family,
and all torus knots are represented by elastic rod centerlines.
\end{abstract}
\maketitle
\section*{Introduction}

The elastic curve of Bernoulli and Euler arises as the solution of the
variational problem of minimizing the total squared (geodesic)
 {\it curvature} (defined below) of a curve
with prescribed boundary conditions and fixed length.
  In \cite{L-S1}, Langer and Singer parametrized the space of
(similarity classes of) elastic curves by a triangular
region of the plane; the closed curves form a countable set of points 
lying on a single curve $c$ within the triangle.  The elastic curves determined
by points on $c$ are {\it quasiperiodic}; that is, they are either
closed curves or they wind densely around a torus of revolution and are
self-congruent under a discrete group of rotations of the torus about its axis.
The closed curves are in 1--to--1 correspondence with the knot types which
are $(k,n)$-torus knots for $k < n/2$.
 
A generalization of the elastic curve can be obtained by adding the
additional constraint that the integral of the {\em torsion} be held
constant.  Thus we seek an extremal $\gamma$ for the Lagrangian
\bel{Lagrange}
{\cal F}[\gamma] = \lambda_1\int_\gamma ds + \lambda_2 \int_\gamma \tau\,ds
+\lambda_3\int_\gamma \k^2 ds,\quad \lambda_3 \ne 0
\end{equation}
Here, $\k$ and $\tau$, the curvature and torsion, are the coefficients 
of the {\it Frenet equations} for a canonical orthonormal (Frenet) frame 
$T,N,B$ along $\gamma$, given by
$$\frac{d\gamma}{ds} = T,\qquad
T' = \k N,\qquad N' = -\k T + \tau B, \qquad B' = -\tau N.$$
The constants $\lambda_1$, $\lambda_2$ and $\lambda_3$ 
are Lagrange multipliers.

That this is an appropriate generalization of elastic curves is a 
consequence of a theorem of Langer and Singer
(cf. \cite{L-S2}, Theorem 1): If $\gamma$ is the centerline of a {\it
uniform symmetric Kirchhoff elastic rod} (see section \ref{sec:stable}),
then $\gamma$ is an extremum for the Lagrangian \eqref{Lagrange}.   
In fact, as we show in section \ref{sec:stable}, the converse of this
theorem also holds.  Thus we may refer to extremals of \eqref{Lagrange}
as elastic rod centerlines.  Our main results are

\begin{thm}\label{knottypethm} Every torus knot type is realized by a smooth 
closed elastic rod centerline.\end{thm}

\begin{thm}\label{homotopythm} The similarity classes of quasiperiodic 
elastic rod centerlines are parametrized by the closed unit disc.  
(Quasiperiodic centerlines are those which are either closed or
wind densely around a torus.)
The closed rod centerlines form a countable family of curves in the disc.
Each such curve represents a regular homotopy of closed curves within
the set of smooth closed elastic rod centerlines.
For any relatively prime positive integers $k,n$,
there exists a regular homotopy, between the $k$-times-covered circle
and the $n$-times covered circle. 
 The homotopy includes exactly one elastic curve, one self-intersecting  
elastic rod centerline, and one closed curve of constant torsion. \end{thm}
\begin{figure}[h]\label{movie}
\begin{center}
\includeps{file=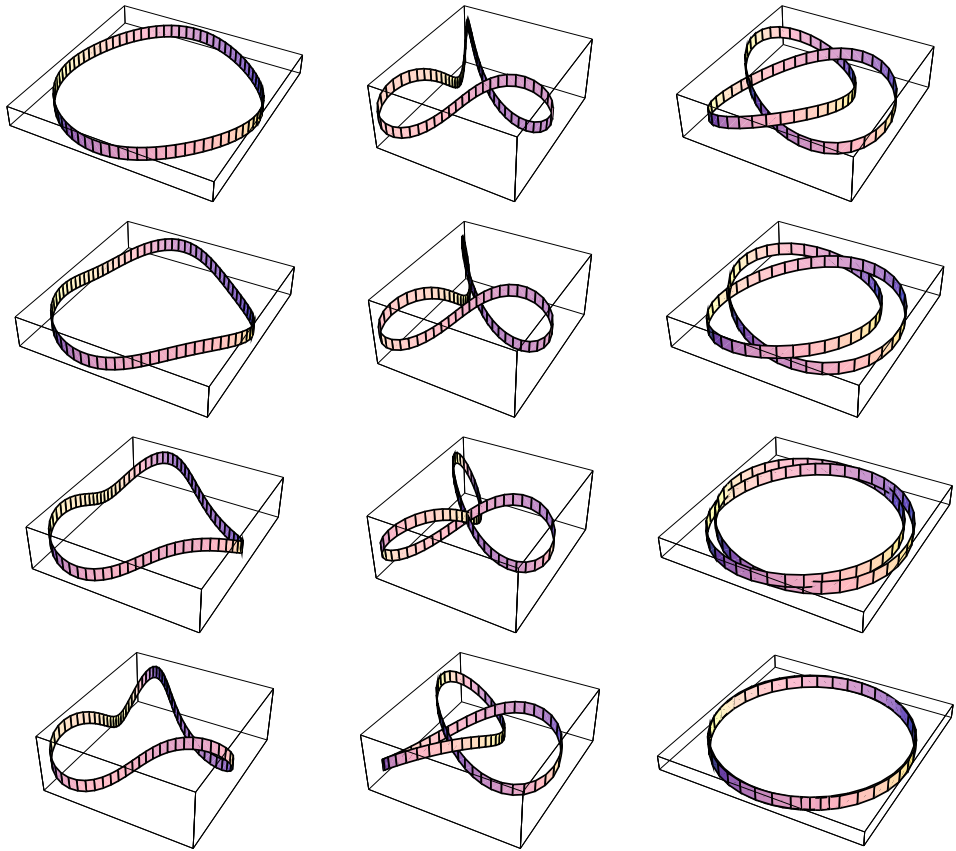, scale=1.3}\end{center}
%%source--notebook "homotopy13" in generalrods
\caption{Elastic rod centerlines along the homotopy joining
the once-covered circle to the twice-covered circle.}
\end{figure}
\noindent
Figure 1 shows the members of one such homotopy of
closed elastic rod centerlines.

We will now give an outline of the rest of the paper.  
In \S1 we detail the explicit solution of the variational problem for
\eqref{Lagrange}.
In \S2 we show how parameters can be chosen so that the equilibria are smooth
closed curves, and we show that, with a few exceptions, these are embedded and
lie on tori of revolution.
In \S3 we introduce a new parametrization of the space $\D$ of quasiperiodic
elastic rods, and point
out loci inside this space that represent elastic curves, curves of constant
torsion, self-intersecting curves, etc.
In \S4 we consider what the elastic rods look like at the boundary of $\D$,
and prove Theorem \ref{knottypethm}.
In \S5 we examine the level curves of a certain function on $\D$ which turn
out to give the homotopies of Theorem \ref{homotopythm}, and prove a
uniqueness theorem 
(Theorem \ref{torsionthm}), conjectured in \cite{C-I}, on the knot
types of elastic rods of constant torsion.
In \S6 we relate the variational problem to the Kirchhoff elastic rod and
discuss the question of stability.  The problem of stability
for closed elastic rods has been extensively investigated, using 
bifurcation theory, by K. Rogers (see \cite{R}, \cite {M-R-M}).
Here we examine the second variation formula in certain specific examples.
This leads to a result in \S7, specifying the precise nature of stability
for the untwisted ``figure eight'' elastic rod.

The authors would like to thank Joel Langer, Annalisa Calini, and
Ron Perline for helpful discussions and encouragement.

\section{Elastic Rods in Cylindrical Coordinates}
In this section, we will recapitulate the integration of the Euler-Lagrange
equations for the
Kirchhoff elastic rod by Langer and Singer \cite{L-S2}.

Suppose $\gamma$ is a space curve which is critical for the functional
\eqref{Lagrange}
with respect to variations that keep the endpoints fixed to second order.  The
resulting Euler-Lagrange equations (see \S4 in \cite{L-S2}) 
have the following first integrals:
\bel{firstint}\begin{aligned}
\k^2(2\lambda_3\tau - \lambda_2) &= c\\
\lambda_3^2\k_s^2 + \tfrac14(\lambda_3\k^2-2\lambda_1)^2
+\k^2(\lambda_3\tau-\lambda_2)^2 &= \mu^2
\end{aligned}\end{equation}
The latter equation allows us to express the curvature as
$$\k^2 = \k_0^2(1-\dfrac{p^2}{w^2}\sn^2(t,p)),\qquad t= \dfrac{\k_0 s}{2w}$$
where $0\le p\le w\le 1$, and $p,w$ are related to the other parameters by
\bel{crel}
\begin{aligned}
\dfrac{4\lambda_1\lambda_3-\lambda_2^2}{\lambda_3^2} &= {\k_0^2 \over w^2}
(3w^2-p^2-1)\\
{c^2 \over \lambda_3^2}&={\k_0^6 \over w^4}(1-w^2)(w^2-p^2)
\end{aligned}\end{equation}
Another consequence of the first integrals is that the vector field
$$J = \tfrac12(\lambda_3\k^2-2\lambda_1)T + \lambda_3 \k_s N 
+ \k(\lambda_3\tau-\lambda_2)B$$
is constant along $\gamma$ (cf. Thm. 2 in \cite{L-S2}).  Setting
$\k=\k_0$ in \eqref{firstint} gives
\bel{Jlength}
|J|^2=\mu^2=\dfrac14(\lambda_3\k_0^2-2\lambda_1)^2
+\dfrac{\k_0^2}4(\dfrac{c}{\k_0^2}-\lambda_2)^2
\end{equation}
In addition to $J$,  
$$I= \lambda_2 T + \lambda_3 \k B$$
is the restriction to $\gamma$ of a Killing field, generating a screw
motion in $\R^3$.  Since $I+J\times \gamma$ is constant, $J$ turns out to 
be parallel to the axis of the screw motion.  In fact, we will translate the
curve so that $I = a J + \gamma \times J$.

The Killing fields $I$ and $J$ define a natural system of cylindrical
coordinates
$r,\theta,z$, in which $J$ points along the positive $z$-axis.  Computing 
$I \cdot J$ in two ways gives
$$a\mu^2  =\tfrac12\lambda_3 c - \lambda_1\lambda_2.$$
Computing $I\cdot I$ in two ways gives
\bel{reqn}
\mu^2(r^2+a^2) = \lambda_2^2 + \lambda_3^2\k^2,\end{equation}
giving us the $r$-coordinate along $\gamma$.  Computing $T \cdot J$ gives
\bel{zderiv}\mu z_s = \tfrac12 \lambda_3\k^2 - \lambda_1\end{equation}
and computing $T \cdot I$ gives
$$\theta_s = \dfrac{a\mu z_s - \lambda_2}{\mu r^2}.$$
Since $r^2$ and $\k^2$ are both of the form $A + B \sn^2 t$, $z(s)$ and
$\theta(s)$ can be computed in terms of elliptic integrals of the second
 and third kind, respectively; we will do so below.

\section{Closure Conditions }
Without loss of generality, we may assume $\lambda_3=1$. By applying a
similarity to a rod we may assume that the maximum curvature $\k_0=1$. 
Then the shape of the 
elastic rods depends on  $\lambda_1$, $\lambda_2$ and the constant
of integration $c$, with $\mu>0$ determined by \eqref{Jlength}.
(We will ultimately use a different set of parameters, though.)
We would like to see
what choices of the parameters lead to closed curves.

Since $\k^2$ is a $2K$-periodic function of $t = k_0 s/(2w)$ 
(where $K=K(p)$ is a complete elliptic integral of the first kind),
equation \eqref{reqn} implies that $r$ is periodic.  
Likewise, \eqref{zderiv} shows that $z_t$ is periodic, and $z$ will
be periodic if and only if $\int_0^{2K} z_t dt=0$; this is equivalent to
\bel{zclosure}w^2-p^2 + w^2\lambda_2^2 = \dfrac{2E(p)}{K(p)} - 1
\end{equation}
(cf. equation (26) in \cite{L-S2}).  So, in order for a rod to have 
periodic $z$ coordinate,
$$A(p) := \dfrac{2E(p)}{K(p)} - 1$$
must not be negative.  Since the behaviour of $A(p)$ as a function of $p$
 will be important in what follows---for example, the sign restriction 
gives us an upper bound on $p$---we will summarize the properties of
 $A(p)$ which we will need.
\begin{prop}\label{Aprop}
For $p \in [0,1]$, $A(p)$ is a smooth monotone decreasing
function of $p$, with $A(0)=1$ and $A(1)=-1$; and,  
$$1-p^2 < A+p^2 < 1$$
for $p \in (0,1)$.\end{prop}
\begin{pf}  The behaviour of $A$ at each end follows from that of $E$ and $K$.
Using differentiation formulas from \cite{B-F}, we have
$$\dfrac{dA}{dp} = 2\left( \dfrac{E-K}{pK} - \dfrac{E}{K^2}\left(
\dfrac{E-(p')^2K}{p(p')^2}\right)\right) = 
-\dfrac2p\left( (1-E/K)^2 + \dfrac{p^2E^2}{(p')^2K^2}\right) <0 $$
and 
$$\dfrac{d}{dp}(A+p^2) = -\dfrac2p\left( p' - \dfrac{E}{Kp'}\right)^2 <0.$$
Finally,
$$2p \dfrac{d}{dp} (A + 2p^2 -1) = 4p^2 - (A + 2p^2 -1)^2/(p')^2.$$
The last equation shows that, after $A + 2p^2 -1=0$ at $p=0$, it remains
positive until $p'=0$.\end{pf}

\def\pmax{p_{\text{max}}}
Let $\pmax \approx .9089085$ be the value of $p$ where $A=0$.  
When $p=\pmax$, the closure condition \eqref{zclosure}
forces $w=p$, hence $c=0$ and the torsion is zero, and $\lambda_2=0$;
then the elastic rod becomes a planar figure-eight elastic curve 
(cf. Figure 1 in \cite{L-S1}.)

When we impose the closure condition \eqref{zclosure}, the formula for the
derivative of $\theta$ becomes
\bel{thetat}\dfrac{d\theta}{dt}=U+\dfrac{N}{1-M \sn^2 t},
\end{equation}
where
$$U=a\mu w,\qquad
 N = \dfrac{{U\over 2}(1-A) - 2\lambda_2 \mu w^3}{p^2+A-U^2}-U,
\qquad M = \dfrac{p^2}{p^2+A-U^2}.$$
(It will follow from calculations in the next section
 that $0 \le A-U^2 < (p')^2$,
so $p^2 < M\le 1$.)  Assuming $M<1$, we 
now use formula 434.01 in \cite{B-F} to compute $\theta$.  The
following identity is necessary for simplifying the result:
\bel{Bident}N^2 = \dfrac{(1-M)(M-p^2)}{M}.\end{equation}
Then
$$\theta(t) = Ut + \sign(N)\left[ \left(\left({E \over
K}-1\right)F(\hat\xi,p') +
E(\hat\xi,p')\right) t  - {1 \over 2i}\ln\dfrac{\Theta_1(t-iF(\hat\xi,p'))}
{\Theta_1(t+iF(\hat\xi,p'))}\right],$$
where 
\bel{defxihat}\hat\xi = \sin^{-1}\left(\sqrt{\dfrac{M-p^2}{(p')^2 M}}\right),
\end{equation}
$p'=\sqrt{1-p^2}$, and $\Theta_1$ is the Jacobi theta function 
\cite{B-F} with period $2K$
and $\Theta_1(K+i K')=0$.
The change in $\theta$ over one period of the other coordinates $r,z$ is
given by
\bel{adeltheta}\Delta\theta = 2KU + 2\sign(N)\left((E-K)F(\hat\xi,p')
 +K E(\hat\xi,p')\right).\end{equation}
The rod will be periodic if the parameters are chosen so that $\Delta\theta$
is a rational multiple of $2\pi$.

\begin{thm}\label{curvethm}An elastic rod centerline of non-constant
curvature can intersect itself only at the origin of
the natural system of cylindrical coordinates.  If this does not occur, and 
the centerline is closed, then it is embedded and lies on an embedded torus
of revolution.
\end{thm}
\begin{pf} Note that the curvature is constant if and only
if $p=0$, giving a planar curve; we will assume $p>0$.  We've observed that
$r^2$ 
and $z$ are $2K$-periodic
functions of $t$.  The minimum value of $r^2$ is non-negative 
(see Lemma 3.1 below) and occurs when $t=0$.

Suppose that $r_{\text{min}} >0$.  Then $r(t)$ is a smooth positive function
with
stationary points at $t=0$ and $t=K$, giving the maximum and minimum values in
each $2K$ period.
Because of the condition \eqref{zclosure},
\bel{rodz}z(t)=\dfrac{1}{\mu w}\dfrac{\Theta'(t)}{\Theta(t)},\end{equation}
where $\Theta$ is the Jacobi theta function with period $2K$ and
$\Theta(i K')=0$.  Then $z(t)$ is a smooth odd
function with stationary points where $\dn^2 t = E/K$, giving the maximum
and minimum values in each period.  Thus, $r(t)$ and $z(t)$ trace out a
simple closed curve in the $rz$ half-plane.  It follows that the rods 
lie on embedded tori of revolution.

We would like to show that the rod is in fact embedded.
Suppose $\theta$ and $\phi$ are longitudinal and latitude on
the embedded torus.  We have already seen that $\phi(t)$ is a monotone
increasing
function along the rod, with $\phi(t+2K) = \phi(t)+2\pi$.
Since $\theta_t$ is periodic, it follows that the rod can intersect itself only
when it closes up smoothly.

In the case where $r_{\min}=0$, we can similarly argue that $r(t)$ and $z(t)$
trace out a closed curve in the $rz$ half-plane, which is smooth except at
the origin.  Replacing $\phi$ with a local coordinate along the smooth part of
this curve, we again see that the rod cannot intersect itself except at
the origin.
\end{pf}

We conclude this section by giving formulae for the rods in Cartesian
coordinates.
Let $\Lambda=\Delta\theta/(2K)$.  Then
\bel{rodxy}\begin{aligned} r\cos\theta &= \sqrt{\dfrac{2Kpp'}{\pi}}\
\dfrac{e^{i\Lambda t}\Theta_1(t\pm i\widehat F) + e^{-i\Lambda
t}\Theta_1(t\mp i\widehat F)}
{2\mu w \Theta(t)H_1(i\widehat F)}\\
r\sin\theta &= \sqrt{\dfrac{2Kpp'}{\pi}}\
\dfrac{e^{i\Lambda t}\Theta_1(t\pm i\widehat F) - e^{-i\Lambda
t}\Theta_1(t\mp i\widehat F)}
{2i\mu w\Theta(t)H_1(i\widehat F)},\end{aligned}\end{equation}
where $\widehat F = F(\hat\xi, p')$, and the upper sign is taken when $N>0$,
the lower sign when $N<0$.

The formula \eqref{rodz} for $z$ bears a striking resemblance to the
representation
of planar elastic curves in terms of theta functions obtained by Mumford
\cite{M}.
In fact, one can check that
as $p \to \pmax$, and our elastic rod becomes planar, \eqref{rodz} and
\eqref{rodxy} 
agree with Mumford's formula in the limit.

\section{The Space $\D$ of Quasiperiodic Elastic Rods}
Because the curvature and torsion are $2K$-periodic in $t$, the curve will
consist of successive congruent segments which join together smoothly; when
we impose the closure condition in the $z$ coordinate, the segments will
be congruent by a rotation in $\theta$.  Thus, the rod centerline will
either close up smoothly after
a finite number of segments, or wind densely around a torus.  We call these
{\it quasiperiodic} elastic rods.

Such curves depend on two parameters; one possible choice of parameters
is the elliptic modulus $p$ and the parameter $\lambda_2$.  
However, we will view the parameter space as
a unit disk in the $XY$ plane, where $X$ and $Y$ are 
related to the other parameters as follows: set
$$X =\lambda_2 w, \qquad Y = \dfrac{cw^2}{\sqrt{1-w^2}},\qquad Z=cw^2.$$
Then, from \eqref{crel} and \eqref{zclosure},
$$w = \sqrt{Y^2 + p^2},\qquad Z = Y\sqrt{(p')^2-Y^2},\qquad X^2 + Y^2 = A(p).$$
(The last equation is equivalent to \eqref{zclosure}.)
If $p \in (0,\pmax)$ then $w>0$, and $X$ and $Y$ are seen to be smooth invertible
functions of $c$ and $p$, tending to the origin as $p$ approaches $\pmax$
and tending to the edge of the unit disk as $p$ approaches zero.
In terms of $X,Y,Z$ and $w$,
the constants used in the rod formulae are:
$$\lambda_1 = {1\over 2} + \dfrac{A-1}{4w^2},\qquad
\mu = \dfrac{\sqrt{(1-A)^2 + 4(Z-wX)^2}}{4w^2}$$
$$U=a\mu w = \dfrac{2w(Z-wX)+X(1-A)}{4\mu w^2}$$ 

We will let $\D$ stand for the open unit disk in the $XY$ plane, 
with the origin omitted for
convenience.  On $\D$, we will also use the angle coordinate $\phi$,
where
$$X=\sqrt{A(p)}\cos\phi,\qquad Y=\sqrt{A(p)}\sin\phi.$$
Thus, $p,\phi$ are also valid coordinates on $\D$.

There are several interesting loci inside $\D$.  First, $X=0$ gives
elastic curves.  Next, $Y=0$ gives elastic rods of constant torsion
$\tau = \lambda_2/2$; these rods are the starting point for the more elaborate
knots of constant torsion produced in \cite{C-I}.  Next, recall that
the Killing field
$$I = \lambda_2 T + \lambda_3 \k B = aJ + \gamma \times J$$
generates the motion of the elastic rod, up to a tangential piece, under
the Localized Induction Equation (cf. \cite{L-S2}, Theorem 3 and Corollary 5.1).
While in general this is a screw motion, it is natural to ask when it is
a pure rotation.  Since $J$ is the translational Killing field, 
defining the $z$ axis for the cylindrical coordinates, $I$ generates a rotation
when $I\cdot J=0$, i.e. when $a=0$.  Setting $a=0$ gives
$$2Zw = X(A-1+2p^2 + 2Y^2),$$
showing that $X$ and $Y$ must have the same sign along this curve.  Substituting
in polar coordinates on the disk gives
$$\sin^2\phi = \dfrac{(A-1+2p^2)^2}{1-A^2 - 2A(A-1+2p^2)}.$$

Kida \cite{K} was the first to investigate space curves that move by a pure
rotation under the LIE,
and for this reason we call this curve in $\D$ the ``Kida curve''.  
It is amusing exercise to verify that, as $p \to 0$, $\sin^2\phi$ approaches
$2/3$ along the Kida curve (see Fig. 2).
Kida also claimed, without proof, that the closed elastic rods
in this family have the knot type of a $(k,n)$ torus knot with $|n/k|>2$.
Our results below confirm this.
\begin{figure}[h]\label{disc}
\begin{center}
\includeps{file=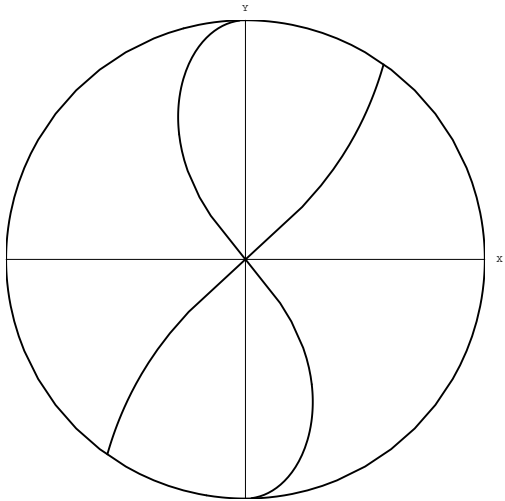, trim=0 1in 0 1in, scale=.4}
%%source notebook "gencircle" in generalrods
\end{center}
\caption{The parameter space $\cal D$ of 
quasi-periodic elastic rods.  Points along the $Y$ axis correspond to elastic
curves, along the $X$ axis to elastic rods of constant torsion, and along the
S-shaped curve to self-intersecting rods.  The Kida curve runs through the
first and third quadrants.}
\end{figure}

A fourth curve in $\D$ of geometric interest corresponds to those rods
which are self-intersecting.  To locate this curve in $\D$, we need
the following
\def\E{V}
\begin{lem} The minimum value of $r^2$ along an elastic rod is
given by 
$$r^2_{\text{min}} =(A - U^2)/(\mu^2 w^2)$$
If $\hat\xi$ is the angle in \eqref{defxihat}, then 
\bel{bigsquare}A - U^2 = (1-p^2)\cos^2\hat\xi = \E^2 \end{equation}
where
$$\E :=(2Xw\sqrt{1-w^2} - Y(1+A-2w^2))/(4\mu w^2).$$
\end{lem}
\begin{pf}Since
$$r^2\mu^2 w^2 = X^2 + w^2 - U^2 - p^2 \sn^2 t,$$
the minimum value of the right-hand side is $A-U^2$. The proof of
\eqref{bigsquare} is a simple (but tedious) calculation using the
equations between various parameters, and we will omit it here.
\end{pf}

 From \eqref{bigsquare}, we see that $r = 0$ is possible only when $\E=0$, i.e. 
$$2Xw\sqrt{1-w^2} = Y(1+A-2p^2-2Y^2)$$
Substituting in the polar coordinates on $\D$ and squaring gives 
\bel{rzero}\cot\phi = \pm\dfrac{A+2p^2-1}{2pp'}.\end{equation}
Substituting back into the equation for $X$ and $Y$ shows
that we should choose the minus sign in \eqref{rzero}, giving
a curve in the second and fourth quadrants of $\D$.  One can check that,
in Figure \ref{disc}, $V$ is positive to the right of this curve and
negative to the
left.

It will be shown in \S 4 that $\E^2 < (p')^2$ on $\D$.  
If we set $\cos\xi = \E/p'$,
then angle $\xi$ is (up to $2\pi$) an analytic function of $X$ and $Y$ that is
equal to $\hat\xi$ to right of the $V=0$ curve and equal to $\pi-\hat\xi$
to the left.  (Remember, $\hat\xi$ only has range $[0,\pi/2]$.)  
It follows from \eqref{Bident} and \eqref{bigsquare} that
$$N^2 = \dfrac{\E^2((p')^2-\E^2)}{\E^2+p^2}.$$
Thus, $N=0$ if and only 
if $\E=0$.  In fact, the sign of $N$ is opposite that of $\E$.  (One can
check this, say, by evaluating $N$ along the Kida curve, where $U=0$.)
Since the term
$$2\left((E-K)F(\hat\xi,p')
 +K E(\hat\xi,p')\right)$$
in \eqref{adeltheta} has limit $\pi$ as $\hat\xi \to \pi/2$, we see that
\bel{bdeltheta}\Delta\theta = 2KU - 2(E-K)F(\xi,p') -2K E(\xi,p')
\end{equation}
\nopagebreak
{\it plus an additional $2\pi$ when $\cos\xi < 0$} (see Figure 3). 
%(See Figure \ref{jagged}
%for a typical plot of
%$\Delta\theta$ as a function of angle $\phi$, for a fixed $p$).
\begin{figure}[h]\label{jagged}
\begin{center}
\includeps{file=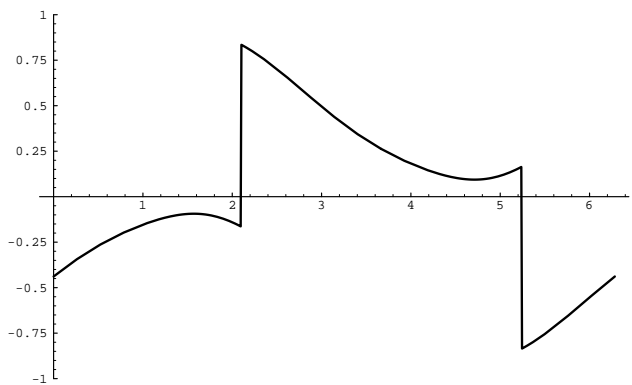, trim=0 2.5in 0 2.5in, scale=.5}
%%source--notebook "homotopy13" in generalrods
\end{center}
\caption{A plot of $\Delta\theta/(2\pi)$ as a function of $\phi$,
for $p=0.8$.  Note the discontinuities, corresponding to rod centerlines
that are self-intersecting.}
\end{figure}

\section{Behaviour on the boundary of $\D$}
\begin{figure}[h]\label{helical}
\includeps{file=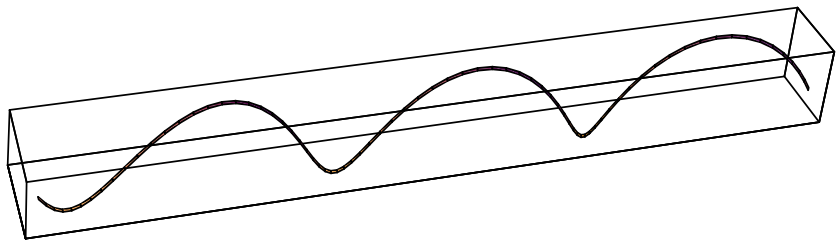, trim=0 3.3in 0 3.3in, scale=.7}
\caption{Three congruent segments of a nearly helical
elastic rod centerline ($p=0.001$, $\phi=\pi/4$).}
\end{figure}

As we approach the boundary of $\D$, the elliptic modulus approaches zero,
and the equation for the curvature
$$\k^2 = 1- \dfrac{p^2}{w^2}\sn^2 t$$
shows that the curvature of the rod becomes constant, provided $w$ has a 
nonzero limit at the boundary point.  (Henceforth we will assume this is
the case; at the points where $p=w=0$, curvature has an undefined limit.)
The torsion has limiting value
$$\tau = (\lambda_2+c/\k_0^2)/(2\lambda_3)=(Z+wX)/w^2$$
But along the boundary, 
$$Z \to Y\sqrt{1-Y^2}, \qquad wX \to X|Y|,$$
so the limiting value of $Z+wX$ is $Y|X|+X|Y|$.  Thus, in the second and fourth
quadrants the torsion goes to zero, and the rod is a circular arc.
In the first and third quadrants the torsion has a nonzero limit, and one might
conclude that the limiting shape of the rod that of a helix.  However, the
rod itself goes off to infinity in the cylindrical coordinates, i.e.
$$r^2_{\text{min}} = \dfrac{A-U^2}{\mu^2 w^2} \to +\infty.$$
One can see this by checking that
$$\lim_{p\to 0} (1-A)/p^2 = 1,\qquad \lim_{p\to 0} 2(Z-wX) = -\cot \phi$$
and hence $U\to 0$ and $\mu \to 0$ along this part of the boundary.
 What one observes experimentally
is that, as the parameters values approach the boundary in the first and third
quadrants, the congruent segments of the rod become shorter (as 
$\Delta\theta$ approaches zero) and more and more helical (see Figure
\ref{helical}).

For what follows, it will be necessary to know how $\Delta\theta$ behaves
along the boundary of $\D$.  It is easy to verify
that, as we approach the part of the boundary that is strictly inside the
second or fourth quadrants,
$$U \to Y,\qquad \mu w^2 \to |XY|,\qquad M\to 0,\qquad N\to Y$$
and so $\theta_t \to 2Y.$
Hence along this part of the boundary, 
\bel{dtedge}\Delta\theta = 2\pi Y.\end{equation}
In particular, as $\phi$ increases, $\Delta\theta$ decreases from $2\pi$ to $0$
along the  
second quadrant edge, and increases from $-2\pi$ to $0$ along the 
fourth quadrant edge.

As we approach the part of the boundary strictly inside the first or third
quadrants,
$$\mu \to 0, \qquad U\to 0, \qquad N \to 0$$
and hence $\theta_t \to 0$.

\begin{pf}[of Theorem \ref{knottypethm}]
Because of \eqref{dtedge} and the continuity of $\Delta\theta$, we see that
for any
rational number $m/n$ such that $|m/n| \in (0,1)$, there are points 
in $\D$, near the
edge in the second ($m/n>0$) or fourth ($m/n <0$) quadrants, where 
$\Delta\theta=2\pi m/n$.  By Theorem \ref{curvethm}, the corresponding
elastic rod will
be embedded and lie on a torus of revolution, and will complete $m$ circuits
around the $z$-axis and $n$ circuits around the waist of the torus before it
closes up smoothly.
\end{pf}

\def\sdt{\widetilde{\Delta\theta}}
\section{Behaviour of $\Delta\theta$ on the interior of $\D$}
Let $\sdt$ stand for the smooth part of 
$\Delta\theta$ (i.e. the right-hand side of \eqref{bdeltheta}).  
In this section we will show that $\sdt$ has
smooth level curves in $\D$ that cross the disk from the second quadrant edge
to the fourth quadrant edge, and which cross each of the $X$-axis, the
$Y$-axis, the
Kida curve (where $U=0$), and the self-intersection locus (where $\E=0$)
 exactly once.
%($\Delta\theta$ jumps by $2\pi$ as we cross that curve.)
We will start by showing that $p$ and $U$ can be used as coordinates
over most of $\D$, and establishing the smoothness of the level curves
in those coordinates.

 From \eqref{defxihat}, we have $A - U^2 = (p')^2 \cos^2\xi$,
and hence for a fixed $p$,
\bel{dudphirel}U\, \di U/\di \phi =(p')^2 \cos\xi \sin\xi\,
\di\xi/\di\phi.\end{equation}

\begin{lem}\label{dudphigood}
For $p \in (0,\pmax)$ fixed, $\di U/\di\phi = 0$
if and only if $\cos\xi=0$.
\end{lem}
\begin{pf} Note that $\E$
is the numerator of $\cos\xi$, and, on the interior $\D$, 
$\E=0$ if and only if $\cos\xi=0$.  Using the derivative formulae
$$\di X/\di\phi = -Y,\qquad \di Y/\di\phi = X, \qquad \di w/\di\phi= XY/w,$$
we obtain
\bel{dudphi}\dfrac{\di U}{\di\phi} = 
\dfrac{\E(2XZ-2wX^2+w(A-1))(2Y(wX-Z)+\sqrt{1-w^2}(A-1))}
{16\mu^2 w^5\sqrt{1-w^2}}.\end{equation}
Note that $\E$
is the numerator of $\cos\xi$, and, on the interior $\D$, 
$\E=0$ if and only if $\cos\xi=0$. 

To complete the proof, we need to check that the other two factors in 
the numerator of $\di U/\di\phi$ cannot be zero.  The middle factor is zero if and only if
$$2XY\sqrt{1-w^2}= w(1+X^2-Y^2);$$
squaring both sides gives
$$w^2(1-A)^2+4p^2X^2=0,$$
which is impossible for $p \in (0,\pmax)$.  The last factor in \eqref{dudphi}
is zero if and only if
$$2wXY= \sqrt{1-w^2}(1-A+2Y^2);$$
squaring both sides gives
$$(1-w^2)(1-A)^2+4Y^2(1-A-p^2)=0,$$
which is also impossible for $p \in (0,\pmax)$, because of Prop. \ref{Aprop}. 
\end{pf}

\begin{cor}\label{xiparam} For $p \in (0,\pmax)$, $\di\cos\xi/\di\phi = 0$
if and only if $U=0$.
\end{cor}
\begin{pf} Since $\cos\xi=\E/p'$, putting \eqref{dudphi} into
\eqref{dudphirel} allows
us to cancel out $\E$ on both sides.  Then the result follows from the
remaining equation.\end{pf}

\begin{cor} For $p \in (0,\pmax)$, $\sin\xi \ne 0$.
\end{cor}
\begin{pf} If $\sin\xi=0$, then by \eqref{dudphirel}, 
either $U=0$ or $\di U/\di \phi=0$.  The
previous lemma rules out the latter, and we just have to verify that, along
the Kida curve, $\cos^2\xi \ne 1$.

When $U=0$, $2wZ = X(A-1+2w^2)$, and hence
$$\cos\xi = \dfrac{A(A-1+2p^2)}{4p'\mu w^2 Y}$$
and
$$16\mu^2w^4 = (1-A)^2 + 4A(1-A-p^2).$$
But squaring each side of $2wZ = X(A-1+2w^2)$ gives
$$A(A-1+2p^2)^2 = Y^2((A-1+2p^2)^2 + 4(A+p^2)(1-A-p^2))=16\mu^2Y^2w^4.$$
Hence $\cos^2\xi = A/(p')^2 <1$.
\end{pf}

By Lemma \ref{dudphigood}, we know that away from the curve $V=0$, $p$ and $U$
can be used as coordinates on $\D$.  Using \eqref{bdeltheta} and 
$(p')^2 \cos^2\xi = A - U^2$,
we compute
\bel{ddu}
\dfrac{1}{2K} \dfrac{\di \sdt}{\di U} = 1 - \dfrac{U(E/K-(p')^2\sin^2\xi)}
{(p')^2 \cos\xi\sin\xi\sqrt{1-(p')^2\sin^2\xi}}.\end{equation}

\begin{lem}\label{dduzero} If $p \in (0,\pmax)$ is fixed and $\cos\xi \ne 0$, 
$\di \sdt/\di U =0$ only when $X=0$.
\end{lem}
\def\w{\omega}
\begin{pf} Let $U=\sqrt{A}\cos\w$ and $\E=p'\cos\xi = \sqrt{A}\sin\w$.
Using $(p')^2 \sin^2\xi = (p')^2 -A\sin^2 \w$, the right-hand side of
\eqref{ddu} is zero if and only if
\bel{dduequiv}\tan\w\sqrt{((p')^2 -A\sin^2 \w)(p^2+A\sin^2\w)} - (E/K-(p')^2
+ A\sin^2\w)=0.
\end{equation}
We will show that the left-hand side of \eqref{dduequiv} is a monotone function of $\w$.
Then it will follow
that the right-hand side of \eqref{ddu}, as a function of angle $\phi$,
can be equal to zero at most once between its vertical asymptotes, 
which occur when $\cos\xi=0$.  

That \eqref{ddu} vanishes when $X=0$ is easy to verify. Let
$$h = p^2+A\sin^2\w.$$
Then, when $X=0$, $Y^2=A$, and 
$$\tan\w = \dfrac{-(1+A-2w^2)}{2w\sqrt{1-w^2}},\qquad h =
\dfrac{(1-A)^2}{16\mu^2w^4},
\qquad 1-h = \dfrac{(1-w^2)(1+A)^2}{16\mu^2 w^4}$$
and
$$16\mu^2 w^4 =(1-A)^2 + 4Y^2(1-w^2).$$
Substituting these in the left hand side of \eqref{dduequiv} gives zero.
 
To complete the proof, let $f=\tan\w\sqrt{h(1-h)}-h$, which
differs from \eqref{dduequiv} by a term depending only on $p$.  We compute
\begin{align*}
\dfrac{\di f}{\di\w} &= \dfrac{h(1-h)\sec^2\w + A(1-2h)\sin^2\w
- 2A\sin\w\cos\w \sqrt{h(1-h)}}{\sqrt{h(1-h)}}\\
 &= \dfrac{(\sqrt{h(1-h)}-A\sin\w\cos\w)^2 +
\tan^2\w(A+p^2)(1-A-p^2)}{\sqrt{h(1-h)}}.
\end{align*}
Since all the terms in the numerator and denominator are nonnegative, 
and never simultaneously zero, we are done.
\end{pf}

\begin{lem}\label{monotone}
For $p\in (0,\pmax)$ fixed, $\sdt$ is a monotone increasing
function of $\phi$ when $X>0$, and monotone decreasing when $X<0$.
\end{lem}
\begin{pf} First, consider just
the $X>0$ side.  Suppose we are at a point where $\cos\xi\ne 0$; then
by (\ref{dudphigood}) and (\ref{dduzero}), $\di\sdt/\di\phi \ne 0$.  On the
other hand,
when $\cos\xi = 0$, this happens in the fourth quadrant, and by (\ref{xiparam}),
$\di\xi/\di\phi \ne 0$ there.  From \eqref{bdeltheta} we have
\begin{align*} \dfrac{1}{2K}\dfrac{\di \sdt}{\di \phi} &= 
\dfrac{\di U}{\di \phi} - \dfrac{\di\xi}{\di\phi}
\left[ \dfrac{(E/K-1)}{\sqrt{1-(p')^2\sin^2\xi}}
+ \sqrt{1-(p')^2\sin^2\xi}\right]\\
&=0 - \dfrac{\di\xi}{\di\phi}\left(\dfrac{A+2p^2-1}{2p}\right) \ne
0.\end{align*}
To get the sign for $\di\sdt/\di\phi$, we check that when $U=0$, 
$\sdt$ is an increasing function of $U$, and, by checking the signs of each
term 
in \eqref{dudphi}, that $U$ is an increasing function of $\phi$ when
$\cos\xi > 0$.

For the $X<0$ side, the argument is similar, except that we find that
$U$ is a decreasing function of $\phi$ when $U=0$.
\end{pf}

\begin{pf}[of Theorem \ref{homotopythm}]
Given $k$ and $n$ coprime, we will show that there exists a curve in $\D$
along which $\sdt = -2\pi k/(k+n)$, running from the fourth quadrant of $\D$
to the second quadrant (see Figure \ref{levels}).  This curve will intersect
the $V=0$ curve once, 
afterwhich $\Delta\theta = 2\pi n/(k+n)$.  At each point along the curve, we
form a closed elastic rod by letting $t$ run through $k+n$ periods of length
$2K$.
(This is the least number at which $\Delta\theta$ is an integer multiple of
$2\pi$.)  
Thus at one end the rods will converge to an $k$-covered circle,
and at the other end to an $n$-covered circle.  

\begin{figure}\label{levels}
\begin{center}
\includeps{file=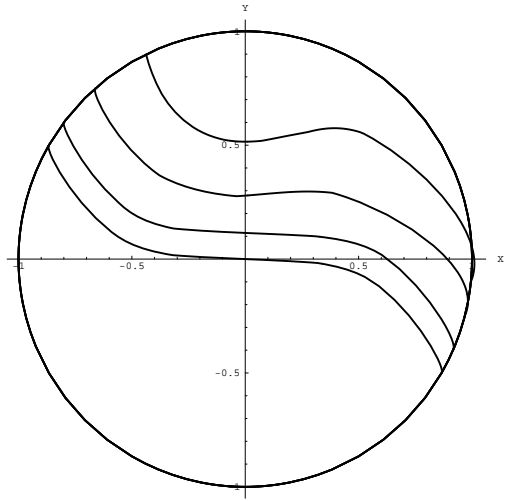, trim=0 1in 0 1in, scale=.4}
%%source notebook "genlevels" in generalrods
\end{center}
\caption{Level curves of $\widetilde{\Delta\theta}$
in $\cal D$.}
\end{figure}

Assume first that $k/(k+n) < 1/2$.
For elastic curves (corresponding to the line $X=0$ in $\D$), it was shown
by Langer and Singer that $\Delta\theta$ is a monotone function of $p$,
decreasing from $0$ to $-\pi$ as we approach the origin along the positive
$Y$-axis (cf. \cite{L-S1}, Figure 2 and Appendix).  
We will follow the level curve along which $\sdt = -2\pi k/(k+n)$ from
a point along the positive $Y$-axis.  Because of 
(\ref{dudphigood}) and (\ref{dduzero}), the level curve is perpendicular to
the $Y$-axis;
however, because of (\ref{monotone}), as we follow it into the first
quadrant, $p$
must be strictly decreasing along the curve.  It cannot veer toward the edge of
the first quadrant, since $\sdt$ has limiting value zero there.  Thus, it
crosses
down into the fourth quadrant, as $p$ continues to decrease.

In the fourth quadrant, the level curve cannot cross the curve $V=0$.  To
see why, note that $V=0$ gives $\xi = \pi/2$ and so
$$\sdt\restr_{V=0} = 2KU - 2\left( (E-K)K' + KE'\right) = 2KU - \pi,$$
using Legendre's relation.  Since $U<0$ in the fourth quadrant, $\sdt <-\pi$ 
along the $V=0$ curve.  Hence the level curve must end at the appropriate point
on the edge of the fourth quadrant.

As we follow this level curve of $\sdt$ in the other direction,
from the positive $Y$-axis into the
second quadrant, we see that $p$ must again be monotone decreasing.  If it 
continues into the third quadrant, it cannot meet the third quadrant edge
because $\Delta\theta=0$ there; if it continues into the fourth
quadrant, $p$ must now be increasing, and the curve is forced into the
origin since it cannot cross the $V=0$ curve.  This too is impossible, since 
$\sdt=-\pi$ at the origin.  Thus, the curve must end at the appropriate
point on the edge of the second quadrant.

It is clear that this level curve must cross the curve $V=0$ at least once
in the second quadrant.  Since \eqref{bigsquare} shows that $U^2=A$
when $V=0$, 
$$\sdt\restr_{V=0} = \pm 2\sqrt{K(2E-K)} - \pi.$$
  Then, since 
$K(2E-K)$ is a monotone decreasing function of $p$, $\sdt$ is
a monotone function along the curve $V=0$, and the level curve crosses $V=0$
only once.

If $k/(k+n) > 1/2$, then we can find a level curve along which $\sdt= -2\pi
n/(k+n)$, passing from the fourth quadrant, through the first, to a point
in the second quadrant where $\Delta\theta = 2\pi -2\pi n/(k+n) = 2\pi
k/(k+n)$.  
Rotating 
this level curve by 180 degrees about the origin, we obtain a level curve of
$\sdt$, passing from a point in the fourth quadrant where $\Delta\theta = 
-2\pi k/(k+n)$, through the third quadrant, to a point in the second quadrant
where $\Delta\theta = 2\pi n/(k+n)$.  (Recall that $\Delta\theta$ is an odd
function on $\D$.)  Once again, this curve intersects the curve $V=0$ exactly
once.

That these homotopies contain exactly one elastic curve and one curve of 
constant torsion follows from the monotonicity of $\Delta\theta$ along the
$X=0$ line, and Theorem \ref{torsionthm} below. 
\end{pf}

The families of elastic rods in $\R^3$ obtained in the proof of Theorem
\ref{homotopythm} all have the
same discrete symmetry group, generated by rotation by $2\pi k/(k+n)$.  Among 
these rods, there is one elastic curve, one rod of constant torsion, one
self-intersecting centerline, and one curve of Kida type moving by a
rotation under the
LIE.  The elliptic modulus is greatest for the elastic curve; since the length,
with respect to the `natural' parameter $t$, is $2(n+m)K$, one might say that
for these rods of fixed symmetry and fixed maximum curvature $\k_0=1$, the
`natural length' is greatest for the elastic curve.  This seems to echo the
mountain-pass argument in \cite{L-S1}, \cite{L-S3},
 which showed that among curves with fixed symmetry and
fixed length, this same elastic curve is the minimax for $\int \k^2 ds$.   

\begin{thm}\label{torsionthm} Given any relatively prime integers $k,n$ such
that $|k/n| <1/2$, there exists a unique smooth closed elastic rod of
constant torsion with the knot type of a $(k,n)$ torus knot.
\end{thm}
\begin{pf} Existence follows from the continuity of $\Delta\theta$ along
the $X$-axis in $\D$, and its limiting values of $-\pi$ as $X \to 0^+$
and $\pi$ as $X \to 0^-$.  To get uniqueness, we need to verify that
$\Delta\theta$ is a monotone function of $p$ along the $X$-axis.  Lemma 
\ref{dudphigood} and \eqref{dudphirel} imply that $p, \xi$ are valid
coordinates 
in the vicinity of this axis.  We will differentiate $\sdt$ with 
respect to $p$, holding $\xi$ fixed.

 From \eqref{bigsquare}, using the differentiation formulas from the proof
of Prop. \ref{Aprop},
$$2U\dfrac{\di U}{\di p} = \dfrac{2}{p}\left( A - \dfrac{E^2}{K^2(p')^2}\right)
+2p \cos^2 \xi.$$
Since $\sdt$ can be written as $\sdt = 2KU - \pi \Lambda_0(\xi, p)$, where
$\Lambda_0$ is Heuman's lambda function (see \cite{B-F}, 150.03), then
\begin{align*}
\dfrac12 \dfrac{\di \sdt}{\di p} &= U\left( \dfrac{E}{p(p')^2}-\dfrac Kp\right)
+\dfrac K{pU}\left( A - \dfrac{E^2}{K^2(p')^2}+p^2\cos^2 \xi \right)
+\dfrac{(K-E)\sin \xi \cos\xi}{p\sqrt{1-(p')^2\sin^2 \xi}}\\
&= \dfrac{(K-E)}{p}\left(\dfrac{(p')^2 \cos^2\xi - 1 - U^2}{2(p')^2 U}
+ \dfrac{\sin \xi \cos\xi}{\sqrt{1-(p')^2\sin^2 \xi}}\right).
\end{align*}
If this is zero, then 
$$(U^2 + 1 - (p')^2 \cos^2\xi)(1-(p')^2 \sin^2\xi) = 2U(p')^2\sin \xi \cos\xi
\sqrt{1-(p')^2\sin^2 \xi},$$
which is equivalent to the vanishing of a sum of squares:
$$\left(U\sqrt{1-(p')^2\sin^2 \xi}-(p')^2\sin \xi \cos\xi\right)^2 + p^2 = 0.$$
Of course, this does not vanish on $\D$.  We have thus shown that, away from
the $U=0$ and $V=0$ curves, $\di\sdt/\di p \neq 0$.  We remark that this
monotonicity calculation includes that of \cite{L-S1} as a special case.
\end{pf}

\section{Elastic Energy and Stability}\label{sec:stable}
\def\Sum{\displaystyle\sum}
In the rest of the paper, we will discuss in more detail the relationship between
critical curves for the geometric functional $\F$ and the centerlines
of elastic rods.  Of particular interest to us will be the question
of which, if any, of the closed nonplanar critical curves for $\F$
are centerlines for {\it stable} elastic rods, i.e., local minima for 
the appropriate energy functional and boundary conditions.

One can model the stress experienced by a rod composed of some
elastic material by associating to the configuration of the rod a {\it material
frame} $(M_1, M_2, T)$.  This is an oriented orthonormal frame such that $T$ is tangent
to the centerline of the rod and $M_1$ and $M_2$ track the twisting of the
material about the centerline.  The components $u_i$ of the Darboux vector
for this frame, defined by
$$T' = -u_2 M_1 + u_1 M_2, \qquad M_1' = u_2 T - u_3 M_2,
\qquad M_2' = -u_1 T + u_3 M_1, $$
are called the {\it material strains}.  (The prime denotes derivative with
respect to arclength.)  In the standard linear model,
 the stress experienced by a uniform, isotropic elastic rod is given
by the energy
\bel{energy}
\EE = \int_0^L \alpha (u_1^2+u_2^2) + \beta u_3^2\, ds =\int_0^L \alpha \k^2 + \beta u_3^2\, ds,
\end{equation}
where $L$ is the length of the rod and $\alpha$ and $\beta$ are 
nonnegative constants depending on the physical characteristics of
the rod.
Equilibria
for $\EE$ are known as {\it Kirchhoff elastic rods}.  When $\beta=0$
as well, $\EE$ only depends on the curvature of the centerline
 of the rod, and equilibrium centerlines are elastic curves.
Equilibrium equations for \eqref{energy} are classical (see \cite{L-M}
for a modern derivation), and imply that 
the twist rate $u_3$ is constant along equilibria.

Now we will formulate a boundary value problem for $\EE$ that
is appropriate to rods with closed centerlines.
Suppose the centerline $\gamma$ is a curve which possesses
a smooth generalized Frenet frame (i.e., satisfying the usual Frenet
equations and extending smoothly across inflection points).
Let $\psi$ be the angle such that $M_1 = N \cos\psi +B\sin\psi$.
Then the total
change $\Delta \psi$ of this angle, measured by integrating
$\psi'$, measures the twisting of 
the material frame relative to the geometry of the centerline.
It is related to the twist rate $u_3$ by the balance equation:
\bel{baal}\int u_3\, ds = \Delta\psi + \int \tau\,ds.\end{equation}
The boundary conditions we will consider are requiring the centerline of
the rod to be closed, and prescribing the value of $\Delta\psi$.
(Physically, this
 can be achieved by attaching collars to the ends of the rod, joining the
ends together, and then adjusting the total twist by rotating
one collar against the other.)  The following theorem clarifies
the relationship between critical rods for this boundary value problem
and critical curves for \eqref{Lagrange}.
%\bel{Lagrange}
%{\cal F}[\gamma] = \lambda_1\int_\gamma ds + \lambda_2 \int_\gamma \tau\,ds
%+\lambda_3\int_\gamma \k^2 ds,\quad \lambda_3 \ne 0
%\end{equation}
\begin{thm}Given any constants $\alpha,\beta>0$, and any curve
$\gamma$ which is critical for a geometric functional of the form \eqref{Lagrange}, there
is a rod with centerline $\gamma$ and twist rate $u_3$ equal to 
a constant $m$,
which is an equilibrium for $\EE$ with respect to variations that preserve a 
certain prescribed value for $\Delta\psi$.  The condition
$\alpha:2\beta m = \lambda_3 : \lambda_2$ determines $m$, 
and the value for $\Delta\psi$ is determined by \eqref{baal}.
\end{thm}
\begin{pf} We will at first assume that $\gamma$ has no inflection points.
Let $L$ be the length of $\gamma$, and let $\cal R$ be the set of rods whose
centerlines are curves of length $L$ with no inflection points, and where
$\Delta\psi$ has some fixed value.  On $\gamma$, construct a rod $\Gamma$ with
this prescribed value and $u_3=m$, determined by
\bel{detm}
L m = \Delta\psi + \int \tau\, ds
\end{equation}

\def\nts{\negthinspace}
Let $\delta$ stand for the gradient operator for functionals on $\cal R$.  
At $\Gamma \in \cal R$,
$$\delta\nts\int (u_3)^2 \, ds= 2m\, \delta\nts\int u_3 \,ds 
= 2m\,\delta\nts\int \tau \, ds,$$ the last equality obtained by 
differentiating \eqref{baal}; hence,
$$\delta\,\EE = \alpha\,\delta\nts\int \k^2\,ds + 2\beta m\,\delta\nts\int\tau\,ds.$$
Since $\gamma$ is a critical curve for $\lambda_2 \int \tau\,ds + \lambda_3
\int \k^2\,ds$ among curves of length $L$, then we see that $\delta\,\EE=0$
if the ratio $\alpha:2\beta m$ is the same as $\lambda_3 : \lambda_2$.

If $\gamma$ has inflection points, $\EE$ and the left hand side of \eqref{baal} are
 still continuous
functions on the space of rods, but $\int\tau\,ds$ can jump by $\pm \pi$ under
smooth variations of the curve (see \cite{M-R} for a simple example).  Thus,
it only makes sense to prescribe $\Delta\psi$ modulo $\pi$.  The gradient
$\delta\int\tau\,ds$ still makes sense, and in fact the above proof goes through
with this simple modification.\end{pf}

We will now outline the computation of the first and second variation
of energy $\EE$.  Under our boundary conditions, the energy reduces
to 
$$\EE = \int \alpha \k^2 + \beta m^2\,ds,$$
where $m$ is the constant determined by \eqref{detm}.

Let $W$ be a periodic vector field along a closed curve $\gamma$,
such that $W$ preserves the length of the curve and respects the
boundary conditions.  Since any length-preserving deformation can
be reparametrized so as to preserve arclength locally along the curve,
we can assume, for the sake of convenience, that $W$ preserves
a unit speed parametrization.
Using variation formulas available
in, for example, \S 6 of \cite{L-S2}, we see that $W'$ must be orthogonal to $T$,
and the $W$-derivative of $m$ is given by
$$m_W = \dfrac{1}{L}\int <W', \k B>\,ds.$$
Then \begin{align*}
\EE_W &= \int <W'', 2\alpha \k N> + 2\beta m<W', \k B>\,ds\\
&= \int<W',-J>\, ds = \int <W, J'> \, ds, \end{align*}
where
$$J = gT +2\alpha \k' N+ (2\alpha \k \tau
-2\beta m\k)B.$$
The tangential component $g$ of $J$ can be chosen so that $J'$ is
normal to the curve and itself defines a length-preserving deformation
of the curve, i.e., $\int <J'', T>\,ds=0$.  This gives 
$$g=\alpha \k^2- \lambda,$$ 
where $\lambda$ is the appropriate constant.  Then, by the usual
reasoning, $J'=0$ for critical curves.

Since the arclength derivative and the $W$-derivative commute, 
$$\EE_{WW} = \int <W, \nabla_W J'>\, ds = \int-<W', \nabla_W J>\,ds.$$
Further computations yield
\begin{align*}
\EE_{WW} &= \int 2\alpha <W'',W''> + 2\beta m( <W',B> <W'',N> - <W'',B><W',N>) \\
&+(\lambda-3\alpha \k^2) <W',W'>\, ds
+ \dfrac{2\beta}{L}\left[ \int <W',\k B>\, ds \right]^2
\end{align*}

Our strategy for investigating stable non-planar rods will be to find 
conditions under which planar elastic rods are unstable.
Letting $W' = \mu N + \nu B$,
and assuming the torsion of the centerline is zero, 
$$\begin{aligned}\label{mnenergy}
\EE_{WW} &= \int 2\alpha(\mu^2\k^2 + (\mu')^2 + (\nu')^2)
+2\be m (\mu'\nu-\nu'\mu)\\
&+(\lambda - 3 \al \k^2)(\mu^2+\nu^2)\, ds+\dfrac{2\beta}{L}\left[ \int_0^L\nu\k\,ds\right]^2.
\end{aligned}$$
Of course, $W'$ must be the derivative of a periodic vector field along
a planar curve, so $\mu$ and $\nu$ are not arbitrary.

\begin{prop}\label{circlestab} An elastic rod with circular centerline is stable
if $|m| < \sqrt{3}\dfrac{2\pi\alpha}{L\beta}$, and unstable if $|m| > 
\sqrt{3}\dfrac{2\pi\alpha}{L\beta}$.
\end{prop}
This result is not new; it was derived, in a different form,
by Zajac \cite{Z}.  For the sake
of completeness, we give our own proof.
\begin{pf} Since curvature is a constant given by $\k = 2\pi/L$,
then $\lambda = \alpha\k^2$.  If we set
$\mu = a_0 + \Sum_{n=1}^{\infty} a_n \cos(nx) + b_n \sin(nx)$ and
$\nu = c_0 + \Sum_{n=1}^{\infty} c_n \cos(nx) + d_n \sin(nx)$, where $x=\k s$, then
the restrictions on $W'$ force $c_0=0$ and $a_1=b_1=0$.  Now,
$$\EE_{WW} =  
2\pi \Sum_{n=2}^{\infty} \left[ \alpha\k(n^2 (a_n^2+b_n^2) + (n^2-1)(c_n^2+d_n^2))
+ 2\beta m n(b_n c_n - a_n d_n) \right]$$
The term in the summation is a sum of two quadratic forms, in $(a_n,d_n)$ and in $(b_n,c_n)$,
which both have the same trace and determinant.  These will be indefinite
if $(\alpha \k)^2 (n^2-1) < \beta^2 m^2$, and definite
otherwise.
\end{pf}

\section{Stability of the Figure Eight}
The other planar closed centerline besides the circle is the figure-eight, where 
$\k = \k_0 \cn(\dfrac{k_0s}{2p},p)$.  The non-constancy of the curvature 
makes it much more difficult to determine under what conditions
$\EE_{WW}$ can be negative.  We will
proceed this way: treating $\EE_{WW}$ as a functional on $L$-periodic functions $\mu(s)$ and
$\nu(s)$ satisfying
\bel{normone}
\int_0^L (\mu^2+\nu^2)\,ds = 1,\end{equation}
we will attempt to minimize $\EE_{WW}$.
The relevant Euler-Lagrange equations are
\bel{varcond}
\begin{aligned} (\lambda-C-\al \k^2)\mu -2\alpha\mu''-2\beta m \nu' &= 0\\
(\lambda-C-3\al \k^2)\nu - 2\al\nu'' +2\beta m \mu' + \dfrac{2\beta \k}{L}
\int_0^L \nu\k \,ds &=0. \end{aligned} \end{equation}
Here, $C$ is a Lagrange multiplier.  One easily computes that, if 
$\mu$ and $\nu$ satisfy these equations and \eqref{normone} as well, then 
$\EE_{WW} = C$.

\def\Lame{Lam\'e}
\def\L{{\cal L}}
\newtheorem{ass}[prop]{Assumption}
The above equations \eqref{varcond} uncouple if we make the
\begin{ass} Our figure-eight elastic rod is untwisted ($m=0$).
\end{ass}

In terms of the independent variable
$t = \k_0 s/(2p)$, \eqref{varcond} becomes 
\begin{align}
\dfrac{d^2 \mu}{dt^2} &= \left(2p^2 \sn^2t - \left(1+
\dfrac{2p^2 C}{\k_0^2\al}\right)\right) \mu \label{mucond}\\
\dfrac{d^2 \nu}{dt^2} &= \left(6p^2\sn^2t - \left(1+4p^2+
\dfrac{2p^2 C}{\k_0^2\al}\right)\right) \nu +
\dfrac{p^2\be}{K\al}\cn t\,\int_0^{4K}\nu\,\cn t\,dt.
\label{nucond}\end{align}
Here, we have used $L=8pK/\k_0$ and 
$\lambda = \k_0^2\al(2p^2-1)/(2p^2)$ for the figure eight.

Recall now \Lame's equation:
\bel{lamey}
\dfrac{d^2 y}{dt^2} = \left(n(n+1)p^2 \sn^2t - h\right)y,
\end{equation}
where $n$ is usually taken to be a positive integer.  The 
eigenvalues---that is, values for $h$ for which 
there exist $4K$-periodic solutions---come in two flavours: 
for the lowest $2n+1$ eigenvalues, there exists
a unique periodic solution (up to multiple), and this is given by a 
{\it \Lame\ polynomial}, a polynomial of degree $n$ in elliptic
functions $\cn,\dn$ and $\sn$.
For the remaining eigenvalues, which will not be of interest to 
us, periodic solutions {\it coexist} \cite{M-W};
that is, there always exist two independent periodic solutions.  (How to 
write down the latter solutions in terms of theta functions
is described in the last few pages of \cite{W-W}.)

\begin{prop}Suppose $W' = \mu N + \nu B$ for a periodic vector field
$W$, and $\mu$ satisfies \eqref{mucond}.  If $C<0$, then $\mu=0$.
If $C=0$ then $\mu$ is a constant multiple of $\k$, and $\EE_{WW}$
depends only on $\nu$.\end{prop}
\begin{pf}\eqref{mucond} is a \Lame\ equation for $n=1$.  According 
to Ince \cite{I}, the solutions of \eqref{lamey}
for the three lowest eigenvalues are
$$\begin{tabular}{|l|l|}\hline
$y=\dn t$ & $h=p^2$ \\ $y=\cn t$ & $h=1$ \\$y=\sn t$ & $h=1+p^2$\\
\hline \end{tabular}$$
If $C<0$, then $h<1$ and $\mu=\dn t$ up to multiple, but this is not possible.
For, if $W=aT + bN + cB$ then 
\bel{abcond}
a'=b\k, \qquad b'=\mu-a\k, \qquad c' = \nu.
\end{equation}
Using $\k=\cn t$, and applying variation of parameters, we find that
$a$ and $b$ will be periodic if and only if 
$\int_0^{4K} \mu\sn t\dn t\,dt=0$ and $\int_0^{4K}\mu (1-2p^2 \sn^2 t)\,dt=0$.  
For $\mu = \dn t$, the former is true but not the latter.

If C=0, $\mu=q \k$ for some constant $q$, and $W-qT$ will be a
variation vector field which has the same action on $\EE$ as $W$,
with $(W-qT)' = \nu B$.\end{pf}

The analysis of the homogeneous equation \eqref{nucond}, for the binormal
component $\nu$ of $W'$, involves the \Lame\ operator for $n=2$:
$$\L:=(d/dt)^2 - (6p^2 \sn^2 t - h).$$
Notice that setting $h=1+4p^2+\dfrac{2p^2 C}{\k_0^2\al}$ makes $\L$
the left-hand side of \eqref{nucond}.
The solutions of $\L y =0$ for the three lowest eigenvalues are:
$$\begin{tabular}{|l|l|}\hline
$y=\sn t \dn t$ & $h=1+4p^2$\\
$y=\cn t \dn t$ & $h=1+p^2$\\
$y=1-(1+p^2-\sqrt{1-p^2+p^4}) \sn^2 t$& $h=2(1+p^2-\sqrt{1-p^2+p^4})$ \\
\hline \end{tabular}$$

Our objective will now be to obtain the value of $\beta/\alpha$ for
which \eqref{nucond} has a nontrivial periodic solution for $C=0$,
and to show that this value is precisely the boundary between 
stability and instability.

\begin{lem} For values of $h$ lying in an open interval that extends to 
$-\infty$ and contains the three lowest eigenvalues (but no other eigenvalues)
of $\L$, there is a unique value of $H$ for which
\bel{Hcond}
\L\nu = H\cn t \int_0^{4K}\nu\cn t\,dt,\end{equation}
has a nonzero periodic solution, and this value $H$ is a continuous, strictly increasing function of $h$.
\end{lem}
\begin{pf} First, assume that $h$ is not an eigenvalue of $\L$.  Then
the equation
\bel{newcond}
\L \nu = \cn t.\end{equation}
has a unique periodic solution, constructed by variation of parameters.
Any solution of \eqref{Hcond} will be a multiple of this, and so $H$
is uniquely determined, in terms of the solution of \eqref{newcond}, by
$$H= \left(\int_0^{4K} \nu \cn t\, dt\right)^{-1}.$$
Furthermore, if $w=\di\nu/\di h$ for the solution 
of \eqref{newcond}, then $\L w = -\nu$.  Using the fact that
$\L$ is self-adjoint,
$$\dfrac{d}{dh}\int_0^{4K} \nu\cn t\, dt = 
 \int_0^{4K} w\,\L \nu\, dt = -\int_0^{4K} \nu\, \L w\, dt=-\int_0^{4K} \nu^2\, dt$$

If $h$ is either the first or third eigenvalue, \eqref{newcond} has
a solution, unique up to addition of an eigenfunction $y$, and $H$
is determined as before.  Moreover, when $h$ is in the vicinity
of either of these eigenvalues, we can construct a family of
solutions to \eqref{newcond} depending continuously on $h$ 
by requiring $v$ to be an odd function
of $t-K$.  (This condition is already met by the non-eigenvalue
solutions.)  Thus, $H$ depends continuously on $h$.

If $h=1+p^2$, no periodic solution to \eqref{newcond} exists.  Instead,
in the vicinity of this eigenvalue, we let $(\nu,H)$ be the
unique solution to $\L \nu = H\cn t$ such that $\nu$ is periodic
and $\int_0^{4K} \nu\cn t\, dt=1$.  Then $H$ depends continuously on
$h$, with $H=0$ when $h=1+p^2$.
\end{pf}

\begin{thm}\label{stabilitythm} The untwisted figure eight elastic rod is stable if 
$\beta/\al > 2$ and unstable if $\beta/\al <2$.
\end{thm}
\begin{pf}  We know, so far, that variations for which $\EE_{WW} <0$
exist only if there is a solution of \eqref{Hcond} for $h <1+4p^2$
and $H=\dfrac{p^2\be}{K\al}$.  Since $H$ is a monotone increasing
function of $h$, this happens only when  
$\dfrac{p^2\be}{K\al}$ is less than the critical value of $H$ 
corresponding to $h=1+4p^2$.  We will now calculate that value.

Let $y_2 = \sn t\dn t$.  Using reduction of order, we obtain a second, linearly
independent solution of the \Lame\ equation $\L y =0$ for $h=1+4p^2$, given by
$$y_1 =\cn^2 t \dn t + \dfrac{p^4}{(p')^2} \cn t\sn^2 t-\left[t+\dfrac{2p^2-1}
{(p')^2}E(t) \right]\dn t\sn t,$$
where $E(t)$ is an incomplete elliptic integral of the second kind
 \cite{B-F}.

 From now on, let the prime denote $d/dt$.  Rewrite the nonhomogenous equation \eqref{newcond} as the system
\bel{linsys}
\dfrac{dz}{dt} = \pmatrix 0 & 1 \\ 6p^2 \sn^2t - h & 0\endpmatrix z 
+ \pmatrix 0\\ \cn t\endpmatrix, \qquad
 z=\pmatrix \nu\\ \nu' \endpmatrix\end{equation}
Let $Y = \pmatrix y_1 & y_2 \\y_1' & y_2' \endpmatrix$; then
variation of parameters gives 
$z=Y\left(\ell + \pmatrix c_1 \\c_2 \endpmatrix\right)$, where
$$\ell(t) = \int_0^t Y^{-1}(x) \pmatrix 0\\ \cn x\endpmatrix dx.$$
To ensure that $z$ is periodic, we must take
$$c_1 = \dfrac{\int_0^{4K} y_1 \cn t\, dt}{y_1'(4K)}=-\tfrac12.$$
(Computations for the numerator and denominator both use the fact that, for the
elliptic modulus $p$ associated to the figure eight elastic rod, the complete
elliptic integrals satisfy $2E = K$.)  
Since the value of $H$ is unaffected by adding multiples of $y_2$
to $\nu$, we can let $c_2=0$. Then
\begin{align*}\int_0^{4K} \nu \cn t\,dt &=\int_0^{4K} ((\ell_1(t) +c_1)
y_1(t)+\ell_2(t)y_2(t))\cn t\, dt \\
& = \int_0^{4K}\int_0^x \cn t \cn x (y_2(t) y_1(x)-y_1(t) y_2(x))\,dx\,dt
+ c_1 \int_0^{4K}y_1(t) \cn t\, dt \\
&= \int_0^{4K} y_1(t)(\sn^2 t +c_1) \cn t\, dt,\end{align*}
using the fact that $\int y_2 \cn t\,dt = \tfrac12 \sn^2 t$.  Finally, we
compute that $\int_0^{4K} \nu \cn t\,dt = K/(2p^2)$, $H = 2p^2/K$, and 
$\beta/\alpha=2$.
\end{pf}

Theorem \ref{stabilitythm} implies that a conjecture made in \cite{L-S1}
and proved in \cite{L-S3},
to the effect that there are no nonplanar stable closed elastic curves,
would be false if extended to elastic rods.  For, the twisted circle,
with $\Delta\psi = 2\pi$ and $m=2\pi/L$, is connected to the untwisted
figure eight rod by a continuous 1-parameter family of rods.  ($\Delta
\psi=2\pi$ throughout the family up until we get to the figure eight,
when $\Delta\psi$ jumps to zero as the inflection points appear.)  
Proposition \ref{circlestab} and Theorem \ref{stabilitythm} together
imply that when $\sqrt{3} < \dfrac{\beta}{\alpha} < 2$, both the
circle and the figure eight rods are unstable for their respective
boundary conditions.  Since these two are the only closed planar elastic
rod centerlines, it follows from the existence of a minimum that
{\it when $\beta/\alpha$ is in the above range
 there exist nonplanar closed elastic rods that are local minima for
the elastic energy subject to the boundary condition $\Delta\psi=2\pi$.}

\bibliographystyle{amsplain}
%\textspaceundo		%comment out to double space

\end{document}